\documentclass[a4paper,12pt]{article}

\usepackage{amsmath, amssymb, amsthm}
\usepackage{mathtools}

\usepackage{hyperref}
\hypersetup{
	colorlinks=true,
	linkcolor=blue,
	urlcolor=cyan,
	citecolor=blue,
}

\usepackage{color}

\definecolor{vio}{rgb}{0.7,0,0.5}

\usepackage{titlesec}
\titleformat{\section}{\bfseries}{\thesection}{1em}{}
\titleformat{\subsection}{\itshape}{\thesubsection}{1em}{}

\newfont{\ctv}{msam10}

\newcommand{\bbox}{\mbox{\ctv \symbol{4}}}
\def\QED{{${}\hfill\bbox$}}
\newenvironment{pf}[1]{\par\vskip1mm{\noindent\it #1.}\ }{\QED\par
\vskip2mm}

\def\bpf{\begin{pf}}
\def\epf{\end{pf}}

\def\expe{\hbox{\rm e}}
\def\ve{\varepsilon}
\def\vp{\varphi}

\def\vrt{\vartheta}
\def\dd{\,\mathrm{d}}
\def\dive{\mathrm{\,div\,}}
\def\sign{\mathrm{\,sign}}
\def\sspan{\mathrm{\,span}}
\def\supess{\mathop{\mathrm{\,sup\,ess}}}
\def\for{\mathrm{\ for\ }}
\def\ale{\mathrm{\ a.\,e.}}

\def\play{\mathfrak{p}}
\def\scal#1{\left\langle #1\right\rangle}
\def\sumj{\sum_{j=1}^\infty}
\def\sumi{\sum_{i=1}^\infty}
\def\sumin{\sum_{i=1}^n}
\def\sumiz{\sum_{i=0}^n}
\def\sumim{\sum_{i=0}^{n-1}}
\def\sumk{\sum_{k=0}^\infty}
\def\oi{^{(i)}}
\def\oj{^{(j)}}
\def\on{^{(n)}}

\def\real{\mathbb{R}}
\def\nat{\mathbb{N}}

\def\nato{\mathbb{N}\cup \{0\}}
\def\HH{\mathcal H}
\def\HN{{\mathcal H}_N}
\def\XX{\mathcal X}
\def\YY{\mathcal Y}
\def\BB{\mathcal B}
\def\hejk{\hat e_{jk}}

\def\io{\int_{\Omega}}
\def\ipo{\int_{\partial\Omega}}

\def\be{\begin{equation}\label}
\def\ee{\end{equation}}

\def\ber{\begin{eqnarray}}
\def\eer{\end{eqnarray}}

\def\bers{\begin{eqnarray*}}
\def\eers{\end{eqnarray*}}

\def\bpf{\begin{pf}}
\def\epf{\end{pf}}

\newtheorem{theorem}{Theorem}[section]
\newtheorem{lemma}[theorem]{Lemma}

\newtheorem{hypothesis}[theorem]{Hypothesis}
\newtheorem{proposition}[theorem]{Proposition}
\newtheorem{definition}[theorem]{Definition}

\begin{document}

\title{Degenerate diffusion in porous media\\ with hysteresis-dependent permeability
\thanks{The support from the European Union’s Horizon Europe research and innovation programme under the Marie Skłodowska-Curie grant agreement No 101102708, from the M\v{S}MT grants 8J23AT008 and 8X23001, and from the GA\v CR project 24-10586S is gratefully acknowledged.}
}

\author{Chiara Gavioli
\thanks{Corresponding author. Faculty of Civil Engineering, Czech Technical University, Th\'akurova 7, CZ-16629 Praha 6, Czech Republic, E-mail: {\tt chiara.gavioli@cvut.cz}.}
\and Pavel Krej\v c\'{\i}
\thanks{Faculty of Civil Engineering, Czech Technical University, Th\'akurova 7, CZ-16629 Praha 6, Czech Republic, E-mail: {\tt Pavel.Krejci@cvut.cz}.}
\thanks{Institute of Mathematics, Czech Academy of Sciences, {\v{Z}}itn{\'a} 25, CZ-11567 Praha 1, Czech Republic, E-mail: {\tt krejci@math.cas.cz}.}
}

\date{}

\maketitle

\begin{abstract}
Hysteresis in the relation between the capillary pressure and the moisture content in unsaturated porous media, which is due to surface tension at the liquid-gas interface, exhibits strong degeneracy in the resulting mass balance equation. Solutions to such degenerate equations have been recently constructed by the method of convexification. We show here that the convexification argument works even if the permeability coefficient depends on the hysteretic saturation. The problem of uniqueness remains open in this case.

\bigskip


\medskip

\end{abstract}


\section{Introduction}\label{sec:intr}

Hysteresis in the pressure-saturation relation in porous media presents mathematical problems of several degrees of difficulty. A basic theory was developed by Visintin in the 1980's, and his main results can be found in the monograph \cite{vis}. Some problems have remained unsolved until recently, namely, degeneracy in the pressure-saturation relation and saturation-dependent permeability. These effects are observed in real experiments and naturally, they have been addressed also by mathematicians. However, there are no published results on hysteresis dependence of the permeability coefficient, unless a rate-dependent correction is included in the hysteresis relation (\cite{bavi1,gosch}), or the hysteresis term is regularized in time (\cite{bavi2}) or in space (\cite[Section~1.9.5]{vis2}).

In this paper, we study the following model problem
\begin{align}\label{e1}
s_t &= \dive \kappa(x,s) \nabla u,\\ \label{e1a}
s &= G[u],
\end{align}
where $u$ represents the pressure, $s$ is the saturation (also called the moisture content) given in terms of a Preisach operator $G$ with initial memory $\lambda$, and $\kappa(x,s)$ is a saturation-dependent permeability. The Preisach operator can reproduce most of the hysteresis effects observed in experiments, as it has been discussed in \cite{colli}. Equation~\eqref{e1} follows from the mass balance and a saturation-dependent variant of the Darcy law. The first to predict that in unsaturated soils the hydraulic conductivity (hence the permeability) would be largely dependent on the soil water content was Buckingham in \cite{buck}. As a result, in the literature this is sometimes referred to as the Darcy-Buckingham law, see e.\,g.\ \cite{GrHa}. Note that Elena El Behi-Gornostaeva in her PhD thesis \cite{gorno} discusses a problem of the type \eqref{e1} with hysteresis-dependent permeability under different (including unilateral) boundary conditions, but with a non-degenerate pressure-saturation relation $s = a_0 u + G[u]$, for some $a_0>0$, replacing \eqref{e1a}.

We study here Problem \eqref{e1}--\eqref{e1a} in a bounded Lipschitzian domain $\Omega \subset \real^N$ and time interval $(0,T)$, with boundary condition
\be{e2}
-\kappa(x,s)\nabla u \cdot n(x) = \gamma(x) (u - u^*) \quad \mbox{on } \ \partial \Omega,
\ee
where $n(x)$ is the unit outward normal vector to $\partial\Omega$ at the point $x \in \partial\Omega$, $\gamma(x) \ge 0$ is the boundary permeability not identically equal to $0$, and $u^* = u^*(x,t)$ is a given boundary pressure. We also prescribe an initial condition $u_0(x)$ such that
\be{e3}
u(x, 0) = u_0(x)\ \ale \mbox{ in } \ \Omega.
\ee
The time evolution described by equation \eqref{e1} is degenerate in the sense that the knowledge of $G[u]_t$ does not give a complete information about $u_t$. In order to control the time derivative of $u$, we have to restrict our considerations to the so-called convexifiable Preisach operators and employ the technique developed in \cite{colli} in the framework of anisotropic Sobolev embedding theorems following \cite{bin}. The main message of this paper is that the convexification technique combined with possibly new anisotropic embedding theorems involving also Orlicz spaces makes it possible to allow the permeability coefficient $\kappa$ to depend on the saturation.

The structure of the present paper is the following. In Section~\ref{stat} we list the definitions of the main concepts including convexifiable Preisach operators, and state the main existence Theorem~\ref{t1}. In Section~\ref{disc} we propose a time discretization scheme with time step $\tau>0$ and derive estimates independent of $\tau$. The crucial convexity estimate \eqref{dp6} is weaker than its counterpart obtained in \cite{colli}, and this is due to the presence of the hysteresis operator in the permeability coefficient. It is therefore necessary to devote Sections~\ref{hilb} and \ref{orli} to a deeper auxiliary material on compact anisotropic embeddings in Sobolev and Orlicz spaces. It is shown in Section~\ref{proo} that this version of compactness is sufficient for passing to the limit as $\tau \to 0$ and for proving that the limit is a solution to the original PDE problem with hysteresis.


\section{Statement of the problem}\label{stat}

We study Problem~\eqref{e1}--\eqref{e2} in variational form
\begin{align}
\io \big(s_t \vrt + \kappa(x,s)\nabla u\cdot\nabla \vrt\big)\dd x + \ipo \gamma(x)(u - u^*)\vrt \dd s(x) &= 0 \ \mbox{ for all }\vrt \in W^{1,2}(\Omega)\cap L^\infty(\Omega), \label{e4} \\
s &= G[u] \ \ale \mbox{ in }\Omega\times(0,T). \label{e5} 
\end{align}

The Preisach operator was originally introduced in \cite{prei}. For our purposes, it is more convenient to use the equivalent variational setting from \cite{book}.

\begin{definition}\label{dpr}
Let $\lambda \in L^\infty(\Omega \times (0,\infty))$ be a given function with the following properties: 
\begin{align}
	\label{ge6b}
	&\exists \Lambda>0 : \ \lambda(x,r) = 0\ \for r\ge \Lambda, \forall x \in \Omega,\\[2mm] \label{ge6}
	&\begin{aligned}
		\exists \bar{\lambda}>0 : \ |\lambda(x_1,r_1) - \lambda(x_2,r_2)| &\le \Big(\bar{\lambda}\,|x_1 - x_2| + |r_1 - r_2|\Big)\ \forall r_1, r_2 \in (0,\infty), \forall x_1,x_2 \in \Omega.
	\end{aligned}
\end{align}
For a given $r>0$, we call the {\em play operator with threshold $r$ and initial memory $\lambda$} the mapping which with a given function $u \in L^1(\Omega; W^{1,1}(0,T))$ associates the solution $\xi^r\in L^1(\Omega; W^{1,1}(0,T))$ of the variational inequality
\be{ge4a}
|u(x,t) - \xi^r(x,t)| \le r, \quad \xi^r_t(x,t)(u(x,t) - \xi^r(x,t) - z) \ge 0 \ \ale \ \forall z \in [-r,r],
\ee
with initial condition
\be{ge5}
\xi^r(x,0) = \lambda(x,r) \ \ale,
\ee
and we denote
\be{ge4}
\xi^r(x,t) = \play_r[\lambda,u](x,t).
\ee
Given a measurable function $\rho :\Omega\times(0,\infty)\times \real \to [0,\infty)$ and a constant $\bar G \in \real$, the Preisach operator $G$ is defined as a mapping $G: L^2(\Omega; W^{1,1}(0,T))\to L^2(\Omega; W^{1,1}(0,T))$ by the formula
\be{ge3}
G[u](x,t) = \bar G + \ \int_0^\infty\int_0^{\xi^r(x,t)} \rho(x,r,v)\dd v\dd r.
\ee
The Preisach operator is said to be {\em regular\/} if the density function $\rho$ of $G$ in \eqref{ge3} belongs to $L^\infty(\Omega\times (0,\infty )\times \real)$, and there exist a constant $\rho_1$ and a decreasing function $\rho_0: \real\to \real$ such that for all $U>0$, all $x,x_1,x_2 \in \Omega$, and a.\,e.\ $(r,v)\in (0,U) \times (-U,U)$ we have
\begin{align}
	&0 < \rho_0(U) < \rho(x,r,v) < \rho_1, \label{ge3a} \\[2mm]
	&|\rho(x_1,r,v) - \rho(x_2,r,v)| \le \bar{\rho}\,|x_1 - x_2|. \label{ge3b}
\end{align}
\end{definition}

In applications, the natural physical condition $s=G[u] \in [0,1]$ is satisfied for each input function $u$ if and only if $\bar G \in (0,1)$ and the additional assumptions
\be{irho}
	\int_0^\infty \int_0^\infty \rho(x,r,v)\dd v\dd r \le 1-\bar G, \quad \int_0^\infty \int_0^\infty \rho(x,r,-v)\dd v\dd r \le \bar G,
\ee
hold for a.\,e.\ $x\in \Omega$. However, the existence result in Theorem~\ref{t1} below holds independently of any specific choice of $\rho$ and, in particular, of \eqref{irho}.

Let us mention the following classical result (see \cite[Proposition~II.3.11]{book}).

\begin{proposition}\label{pc1}
	Let $G$ be a regular Preisach operator in the sense of Definition~\ref{dpr}. Then for every $p \in [1,\infty)$ it can be extended to a Lipschitz continuous mapping $G: L^p(\Omega; C[0,T]) \to L^p(\Omega; C[0,T])$.
\end{proposition}

The Preisach operator is rate-independent. Hence, for input functions $u(x,t)$ which are monotone in a time interval $t\in (a,b)$, a regular Preisach operator $G$ can be represented by a superposition operator $G[u](x,t) = B(x, u(x,t))$ with an increasing function $u \mapsto B(x, u)$ called a {\em Preisach branch\/}. Indeed, the branches may be different at different points $x$ and different intervals $(a,b)$. The branches corresponding to increasing inputs are said to be {\em ascending\/} (the so-called wetting curves in the context of porous media), the branches corresponding to decreasing inputs are said to be {\em descending\/} (drying curves).

\begin{definition}\label{dpc}
	Let $U>0$ be given. A Preisach operator is said to be {\em uniformly counterclockwise convex on $[-U,U]$\/} if for all inputs $u$ such that $|u(x,t)|\le U$ a.\,e., all ascending branches are uniformly convex and all descending branches are uniformly concave.
	
	A regular Preisach operator $G$ is called {\em convexifiable\/} if for every $U>0$ there exist a uniformly counterclockwise convex Preisach operator $P$ on $[-U,U]$, positive constants $g_*(U),g^*(U),\bar{g}(U)$, and a twice continuously differentiable mapping $g:[-U,U] \to [-U,U]$ such that
\be{hg}
	g(0)=0, \quad 0<g_*(U) \le g'(u) \le g^*(U), \quad |g''(u)| \le \bar g(U)\ \ \forall u\in [-U,U],
\ee
	and $G = P\circ g$. Here $g'$ and $g''$ denote the first and second derivative of the function $g$.
\end{definition}

A typical example of a uniformly counterclockwise convex operator is the so-called {\em Prandtl-Ishlinskii operator\/} characterized by positive density functions $\rho(x,r)$ independent of $v$, see \cite[Section~4.2]{book}. Operators of the form $P\circ g$ with a Prandtl-Ishlinskii operator $P$ and an increasing function $g$ are often used in control engineering because of their explicit inversion formulas, see \cite{al,viso,kk}. They are called the {\em generalized Prandtl-Ishlinskii operators\/} (GPI) and represent an important subclass of Preisach operators. Note also that for every Preisach operator $P$ and every Lipschitz continuous increasing function $g$, the superposition operator $G = P\circ g$ is also a Preisach operator, and there exists an explicit formula for its density, see \cite[Proposition~2.3]{error}. Another class of convexifiable Preisach operators is shown in \cite[Proposition~1.3]{colli}.

The hypotheses on $\gamma,u^*,\kappa$ in \eqref{e4} are only technical and can be stated as follows.

\begin{hypothesis}\label{hy2}
The boundary permeability $\gamma$ belongs to $L^\infty(\partial\Omega)$, and is such that $\gamma(x)\ge 0$ a.\,e.\ and $\ipo \gamma(x)\dd s(x) > 0$. The boundary source $u^*$ belongs to $L^\infty(\partial \Omega\times(0,T))$, in particular, there exists $U^*>0$ such that $|u^*(x,t)| \le U^*$ a.\,e.; additionally, $u^*_t \in L^2(\partial \Omega\times(0,T))$. The permeability $\kappa:\Omega \times \real \to \real$ is a bounded Lipschitz continuous function, more precisely, there exist constants $\kappa_*, \kappa^*,\bar{\kappa}$ such that for all $s, s_1, s_2\in \real$ and all $x, x_1, x_2 \in \Omega$ we have
\be{hka}
0< \kappa_* \le \kappa(x,s) \le \kappa^*, \quad |\kappa(x_1,s_1) - \kappa(x_2,s_2)| \le \bar\kappa\big(|x_1 - x_2| + |s_1 - s_2|\big).
\ee
\end{hypothesis}

Note that even a local solution to Problem \eqref{e1}--\eqref{e1a} may fail to exist if for example $\lambda(x,r) \equiv 0$ and $\!\dive \big(\kappa(x,s_0)\nabla u_0\big)\ne 0$, and we need an initial memory compatibility condition which we state here following \cite{colli}. A more detailed discussion on this issue can be found in the introduction to \cite{colli}.

\begin{hypothesis}\label{hy1}
Let the initial memory $\lambda$ and the Preisach density function $\rho$ be as in Definition~\ref{dpr}. The initial pressure $u_0$ belongs to $W^{2,\infty}(\Omega)$, and there exist a constant $L>0$ and a function $r_0 \in L^\infty(\Omega)$ such that, for $\Lambda>0$ as in \eqref{ge6b}, $\supess_{x\in \Omega}|u_0(x)| \le \Lambda$ and the following initial compatibility conditions hold:
	\begin{align} \label{c0}
		\lambda(x,0) &= u_0(x) \ \ale \textup{ in } \Omega,\\ \label{c0a}
		s_0(x) &= G[u](x,0) = \bar G + \ \int_0^\infty\int_0^{\lambda(x,r)} \rho(x,r,v)\dd v\dd r \ \ale \textup{ in } \Omega,\\ \label{c1}
		\frac1L \sqrt{\big|\dive \big(\kappa(x,s_0)\nabla u_0\big)}\big| &\le r_0(x) \le \Lambda \ \ale \textup{ in } \Omega, \\ \label{c2}
		-\frac{\partial}{\partial r} \lambda(x,r) &\in \sign\Big(\!\dive \big(\kappa(x,s_0)\nabla u_0\big)\Big) \ \ale \textup{ in } \Omega\ \for r\in (0,r_0(x)), \\ \label{c2a}
		-\kappa(x,s_0(x))\nabla u_0 \cdot n(x) &= \gamma(x) (u_0(x) - u^*(x,0)) \ \ale\ \mbox{\em on } \partial\Omega.
	\end{align}
\end{hypothesis}

Unlike \cite{colli}, here we do not need to assume $\!\dive\!\big(\kappa(x,s_0(x))\nabla u_0(x)\big) \in L^\infty(\Omega)$ since it follows from the fact that $u_0 \in W^{2,\infty}(\Omega)$ together with assumptions \eqref{ge6}, \eqref{ge3b}, and \eqref{hka}. Our main existence result reads as follows.

\begin{theorem}\label{t1}
Let Hypotheses~\ref{hy2} and \ref{hy1} hold, and let $G$ be a convexifiable Preisach operator	in the sense of Definition~\ref{dpc}. Then there exists a solution $u \in L^\infty(\Omega\times (0,T))$ to Problem~\eqref{e4}--\eqref{e5}, \eqref{e3} such that $\nabla u \in L^2(\Omega\times (0,T);\real^N)$, and both $u_t$ and $s_t = G[u]_t$ belong to the Orlicz space $L^{\Phi_{log}}(\Omega\times (0,T))$ generated by the function $\Phi_{log}(v) = v\log(1+v)$.
\end{theorem}

Basic properties of Orlicz spaces are summarized below in Section~\ref{orli}.


\section{Time discretization}\label{disc}

We proceed as in \cite{colli}, choose a discretization parameter $n \in \nat$, define the time step $\tau = T/n$, and replace \eqref{e4}--\eqref{e5} with its time discrete system for the unknowns $\{u_i: i = 1, \dots, n\} \subset W^{1,2}(\Omega)$ of the form
\be{dis1}
\io \left(\frac1\tau(G[u]_i - G[u]_{i-1})\vrt + \kappa(x,G[u]_i)\nabla u_i\cdot\nabla\vrt\right)\dd x + \ipo \gamma(x)(u_i - u^*_i)\vrt \dd s(x) = 0
\ee
for every test function $\vrt \in W^{1,2}(\Omega)$, where $u^*_i(x) = u^*(x,i\tau)$ for $i\in \{1,\dots,n\}$. Here, the time-discrete Preisach operator $G[u]_i$ is defined for an input sequence $\{u_i : i\in\nat\cup\{0\}\}$ by a formula of the form \eqref{ge3}, namely,
\be{de3}
G[u]_i(x) = \bar G + \ \int_0^{\infty}\int_0^{\xi^r_i(x)} \rho(x,r,v)\dd v\dd r,
\ee
where $\xi^r_i$ denotes the output of the time-discrete play operator
\be{de4}
\xi^r_i(x) = \play_r[\lambda,u]_i(x)
\ee
defined as the solution operator of the variational inequality
\be{de4a}
|u_i(x) - \xi^r_i(x)| \le r, \quad (\xi^r_i(x) - \xi^r_{i-1}(x))(u_i(x) - \xi^r_i(x) - z) \ge 0 \quad \forall i\in \nat \ \ \forall z \in [-r,r],
\ee
with a given initial condition
\be{de5}
\xi^r_0(x) = \lambda(x,r) \ \ale
\ee
similarly as in \eqref{ge4a}--\eqref{ge5}. Note that the discrete variational inequality \eqref{de4a} can be interpreted as weak formulation of \eqref{ge4a} for piecewise constant inputs in terms of the Kurzweil integral, and details can be found in \cite[Section 2]{ele}.

For each $i\in \{1,\dots,n\}$, there is no hysteresis in the passage from $u_{i-1}$ to $u_i$, so that \eqref{dis1} is a standard  quasilinear elliptic equation. The existence of a solution follows from a classical argument based on Fourier expansion into eigenfunctions of the Laplacian, Brouwer degree theory, and a homotopy argument. An introduction to topological methods for solving nonlinear partial differential equations can be found in \cite[Chapter~V]{fuku}.


\subsection{Uniform upper bounds}\label{unif}

In this and the following subsections, unless otherwise specified, the letter $C$ stands for a positive constant that may vary from line to line. Following \cite{hilp}, the first idea is to test \eqref{dis1} by $\vrt =H_\ve(u_i - U)$, with $H_\ve$ being a Lipschitz regularization
\begin{equation*}
	H_\ve (s) = \left\{
	\begin{array}{ll}
		0 & \for s \le 0,\\
		\frac{s}{\ve} & \for s\in (0,\ve),\\
		1 & \for s \ge \ve,
	\end{array}
	\right.
\end{equation*}
of the Heaviside function for some $\ve > 0$, and $U = \max\{U^*,\Lambda\}$ with the notation from Hypotheses~\ref{hy2} and \ref{hy1}. We then let $\ve$ tend to $0$. The elliptic term and the boundary term give a non-negative contribution, and we get for all $i \in \{1,\dots,n\}$ that
\be{hil}
\io (G[u]_i - G[u]_{i-1})H(u_i - U)\dd x \le 0.
\ee
We define the functions
\be{psi}
\psi(x,r,\xi) \coloneqq \int_0^\xi\rho(x,r,v)\dd v, \quad \Psi (x,r,\xi) \coloneqq \int_0^\xi v \rho(x,r,v)\dd v.
\ee
In terms of the sequence $\xi^r_i(x) = \play_r[\lambda,u]_i$ we have
\be{Gi}
G[u]_i(x) = \bar G + \int_0^\infty \psi(x,r,\xi^r_i(x))\dd r.
\ee
Choosing $z= U - (U-r)^+ = \min\{U,r\}$ in \eqref{de4a} and using the fact that $\psi$ is an increasing function of $\xi$, we get for all $i\in\nat$, all $r>0$, and a.\,e.\ $x \in \Omega$ that
$$
\big(\psi(x,r,\xi^r_i(x)) - \psi(x,r,\xi^r_{i-1}(x))\big)\big((u_i(x)-U) - (\xi^r_i(x)-(U-r)^+)\big) \ge 0.
$$
The Heaviside function $H$ is nondecreasing, hence,
$$
\big(\psi(x,r,\xi^r_i(x))-\psi(x,r,\xi^r_{i-1}(x))\big)\big(H(u_i(x) - U) - H(\xi^r_i(x) - (U-r)^+)\big) \ge 0.
$$
From \eqref{hil} it follows that
\be{hil2}
\io\int_0^\infty\big(\psi(x,r,\xi^r_i(x)) - \psi(x,r,\xi^r_{i-1}(x))\big) H(\xi^r_i(x) - (U-r)^+)\dd r \dd x \le 0.
\ee
We now proceed by induction over $i=1,\dots n$ and assume that
\be{indu1}
\xi^r_{i-1}(x) \le (U - r)^+ \ \ale
\ee
for some $i\in \{1,\dots,n\}$. We have by \eqref{ge6b}--\eqref{ge6} that $\xi^r_0(x) = \lambda(x,r) \le (\Lambda - r)^+ \le (U-r)^+$ a.\,e., hence \eqref{indu1} is satisfied for $i=1$. If \eqref{indu1} holds for some $i\ge 1$, we get from \eqref{hil2} that
\be{hil3}
\io\int_0^\infty\big(\psi(x,r,\xi^r_i(x)) {-} \psi(x,r,\xi^r_{i-1}(x))\big) \big(H(\xi^r_i(x) {-} (U{-}r)^+) {-} H(\xi^r_{i-1}(x) - (U{-}r)^+)\big)\dd r \dd x \le 0.
\ee
By monotonicity of both functions $\psi$ and $H$, the expression under the integral sign in \eqref{hil3} is non-negative almost everywhere, hence it vanishes almost everywhere, and we conclude that $\xi^r_i(x) \le (U - r)^+$ a.\,e.\ for all $r\ge 0$ and $i\in \{0, 1,\dots,n\}$. Similarly, replacing $H(u_i - U)$ with $H^*(u_i +U)$ in \eqref{hil}, where $H^*(v) = -H(-v)$ for $v \in \real$ and putting $z= (U-r)^+ - U = -\min\{U,r\}$ in \eqref{de4a} we get $\xi^r_i(x) \ge -(U - r)^+$, so that
\be{Lam}
|u_i(x)| \le U, \quad |\xi^r_i(x)| \le (U - r)^+
\ee
for a.\,e.\ $x \in \Omega$ and all $r\ge 0$ and $i\in \{0, 1,\dots,n\}$. In particular, this implies that $\xi^r_i(x) = 0$ a.\,e.\ for all $r\ge U$, and in \eqref{de3} we can actually write
$$
G[u]_i(x) = \bar G + \int_0^U\int_0^{\xi^r_i(x)} \rho(x,r,v) \dd v\dd r,
$$
so that, even if we do not assume the a priori boundedness of $G$ as in \eqref{irho}, by assumption \eqref{ge3a} we obtain
\begin{equation}\label{GiL}
	 \quad |G[u]_i| \le C
\end{equation}
with a constant $C>0$ independent of $i$ and $\tau$.

We further test \eqref{dis1} by $\vrt = u_i$ and get for $s_i = G[u]_i$
\be{des2}
\frac1{\tau}\io (G[u]_i - G[u]_{i-1})u_i\dd x + \io \kappa(x,s_i)|\nabla u_i|^2\dd x + \ipo \gamma(x)(u_i - u^*_i)u_i \dd s(x) = 0
\ee
for all $i\in\{1,\dots,n\}$. Choosing in \eqref{de4a} $z = 0$ and using the fact that the function $\psi$ in \eqref{psi} is increasing, we obtain in both cases $\xi^r_i\ge \xi^r_{i-1}$ or $\xi^r_i\le \xi^r_{i-1}$ the inequalities
$$
(\psi(x,r,\xi^r_i)-\psi(x,r,\xi^r_{i-1})) u_i \ge (\psi(x,r,\xi^r_i)-\psi(x,r,\xi^r_{i-1}))\xi^r_i \ge \Psi(x,r,\xi^r_i)-\Psi(x,r,\xi^r_{i-1}).
$$
Identity \eqref{des2} together with Hypothesis~\ref{hy2} then yield
\be{energ}
\frac1\tau \io\int_0^U \big(\Psi(x,r,\xi^r_i)-\Psi(x,r,\xi^r_{i-1})\big)\dd r \dd x +\io |\nabla u_i|^2\dd x + \ipo \gamma(x)|u_i|^2 \dd s(x) \le C
\ee
for $i \in \{1,\dots,n\}$, with a constant $C>0$ independent of $\tau$. We now sum up inequality \eqref{energ} over $i$. On the one hand, it holds that $\Psi(x,r,\xi^r_n(x)) \ge 0$ a.\,e.; on the other hand, an estimate for the initial time term
$$
\int_0^U \Psi(x,r,\xi^r_0(x)) \dd r = \int_0^U\int_0^{\lambda(x,r)} v\rho(x,r,v) \dd v\dd r \le \frac{\rho_1}{2} \int_0^\Lambda \lambda^2(x,r) \dd r \le C
$$
follows from the definition of $\xi^r_0$ in \eqref{de5} and the assumptions on $\rho$ and $\lambda$ in Definition~\ref{dpr}. Finally, this yields the energy estimate
\be{energy}
\tau\sumiz \left(\io |\nabla u_i|^2\dd x + \ipo \gamma(x)|u_i|^2 \dd s(x)\right) \le C
\ee
with a constant $C>0$ independent of $\tau$.


\subsection{Convexity estimate}\label{conv}

Recall that the operator $G$ is convexifiable in the sense of Definition~\ref{dpc}, that is, for every $U>0$ there exists a twice continuously differentiable mapping $g:[-U, U] \to [-U, U]$ such that $g(0) = 0$, $0 < g_* \le g'(u) \le g^* < \infty$, $|g''(U)| \le \bar{g}$, and $G$ is of the form
\be{ne0}
G = P \circ g,
\ee
where $P$ is a uniformly counterclockwise convex Preisach operator on $[-U,U]$.

Let us fix $U$ from \eqref{Lam} and the corresponding function $g$. The following result is a straightforward consequence of \cite[Proposition 3.6]{colli}.

\begin{proposition}\label{pc}
Let $P$ be a uniformly counterclockwise convex Preisach operator on $[-U,U]$ whose discrete version is as in \eqref{de3}, and let $f$ be an odd increasing function such that $f(0) = 0$. Then there exists $\beta>0$ such that for every sequence $\{w_i: i = -1, 0, \dots, n-1\}$ in $[-U,U]$ we have
\be{Pi}
\begin{aligned}
&\sumim (P[w]_{i+1} - 2P[w]_i + P[w]_{i-1})f(w_{i+1} - w_i) + \frac{P[w]_0 - P[w]_{-1}}{w_0 - w_{-1}}F(w_0 - w_{-1}) \\[2mm]
&\qquad \ge \frac{\beta}{2}\sumim \Gamma(w_{i+1} - w_i),
\end{aligned}
\ee
where we set for $w\in \real$
\be{gf1}
F(w) \coloneqq \int_0^w f(v)\dd v, \qquad \Gamma(w) \coloneqq |w|(wf(w) - F(w)) = |w|\int_0^{|w|}vf'(v)\dd v.
\ee
\end{proposition}

We need to define a backward step $u_{-1}$ satisfying the strong formulation of \eqref{dis1} for $i=0$, that is,
\be{e7}
\frac1\tau(G[u]_0(x) - G[u]_{-1}(x)) = \dive \big(\kappa(x,s_0)\nabla u_0\big) \ \mbox{ in } \Omega
\ee
with boundary condition \eqref{c2a}.
Repeating the argument of \cite[Proposition~3.3]{colli} we use assumptions \eqref{ge3a} and \eqref{c2} to find for each $0<\tau <\rho_0(U)/2L^2$ functions $u_{-1}$ and $G[u]_{-1}$ satisfying \eqref{e7} as well as, thanks to \eqref{c1}, the estimate
\be{inim}
\frac1\tau |u_0(x) - u_{-1}(x)| \le C
\ee
with a constant $C>0$ independent of $\tau$ and $x$. The discrete equation \eqref{dis1} extended to $i=0$ has the form
\be{dp1}
\io \left(\frac1\tau(P[w]_i - P[w]_{i-1})\vrt + \kappa(x,s_i)\nabla u_i\cdot\nabla\vrt\right)\dd x + \ipo \gamma(x)(u_i - u^*_i)\vrt \dd s(x) = 0
\ee
with $w_i = g(u_i)$, $s_i = G[u]_i$ for $i\in \{0,1,\dots,n\}$ and for an arbitrary test function $\vrt \in W^{1,2}(\Omega)$. We proceed as in \cite{colli} and test the difference of \eqref{dp1} taken at discrete times $i+1$ and $i$
\begin{equation}\label{dp2}
\begin{aligned} 
&\io \left(\frac1\tau\big(P[w]_{i+1} - 2P[w]_i + P[w]_{i-1}\big)\vrt + \Big(\kappa(x,s_{i+1})\nabla u_{i+1}-\kappa(x,s_i) \nabla u_i\Big) \cdot\nabla\vrt\right)\dd x\\
&\qquad + \ipo \gamma(x)(u_{i+1} -u_i)\vrt \dd s(x) = \ipo \gamma(x)(u^*_{i+1} - u^*_i)\vrt \dd s(x)
\end{aligned}
\end{equation}
by $\vrt = f(w_{i+1} - w_i)$ with 
\be{dp2f}
f(w) := \frac{w}{\tau + |w|}.
\ee
In agreement with \eqref{gf1} we have
\begin{align}\label{ff}
F(w) &= |w| - \tau\log\left(1 + \frac{|w|}{\tau}\right), \\ \label{fgamma}
\Gamma(w) &= \tau |w|\left(\log\left(1 + \frac{|w|}{\tau}\right) - \frac{|w|}{\tau + |w|}\right).
\end{align}
The hysteresis term in \eqref{dp2} is estimated from below by virtue of Proposition~\ref{pc} as follows:
\begin{equation}\label{dp3}
\begin{aligned}
\frac1\tau \sumim &(P[w]_{i+1} - 2P[w]_i + P[w]_{i-1})f(w_{i+1} - w_i) + \frac{1}{\tau}\frac{P[w]_0-P[w]_{-1}}{w_0-w_{-1}}F(w_0-w_{-1}) \\
&\ge \frac{\beta}{2\tau}\sumim \Gamma(w_{i+1} - w_i).
\end{aligned}
\end{equation}
Note that for $|w| \ge \tau(\expe^2 - 1)$ we have
$$
\frac{|w|}{\tau + |w|} < 1 \le \frac12\log\left(1 + \frac{|w|}{\tau}\right),
$$
so that
$$
\Gamma(w) \ge \frac{\tau}{2}|w|\log\left(1 + \frac{|w|}{\tau}\right).
$$
We denote by $J$ the set of all $i \in \{0,1,\dots, n-1\}$ such that $|w_{i+1} - w_i| \ge \tau(\expe^2 - 1)$, and by $J^\perp \coloneqq \{0,1,\dots, n-1\} \setminus J$ its complement. We thus have
\begin{align}\label{dp3a}
\frac12\sum_{i\in J} |w_{i+1} - w_i|\log\left(1 + \frac{|w_{i+1} {-} w_i|}{\tau}\right) &\le \frac1\tau\sum_{i\in J}\Gamma(w_{i+1} - w_i), \\ \label{dp3c}
\frac12\sum_{i\in J^\perp} |w_{i+1} - w_i|\log\left(1 + \frac{|w_{i+1} {-} w_i|}{\tau}\right) &\le T(\expe^2 - 1),
\end{align}
hence,
$$
\frac{\beta}{2\tau}\sumim \Gamma(w_{i+1} - w_i) \ge \frac\beta{4}\sumim |w_{i+1} - w_i|\log\left(1 + \frac{|w_{i+1} {-} w_i|}{\tau}\right) - C
$$
with a constant $C>0$ independent of $x$ and $\tau$. Moreover, for $w_0 \neq w_{-1}$ we have
$$
0 < \frac{F(w_0-w_{-1})}{|w_0-w_{-1}|} \le 1,
$$
which yields, together with identity \eqref{e7} and assumption \eqref{c1},
\be{dp3b}
0 < \frac{1}{\tau}\left|\frac{P[w]_0-P[w]_{-1}}{w_0-w_{-1}}\right| F(w_0-w_{-1}) \le C
\ee
with a constant $C>0$ independent of $x$ and $\tau$. For $w_0 = w_{-1}$, we interpret $(P[w]_0-P[w]_{-1})/(w_0-w_{-1})$ as $B'_+(w_{-1})$ or $B'_-(w_{-1})$, according to the notation for the Preisach branches introduced after Proposition~\ref{pc1}, see \cite{colli} for more details.
From \eqref{dp3}--\eqref{dp3b} we thus get
\be{dp4}
\frac1\tau \sumim (P[w]_{i+1} - 2P[w]_i + P[w]_{i-1})f(w_{i+1} {-} w_i)
\ge \frac{\beta}{4}\sumim |w_{i+1} {-} w_i|\log\left(1 + \frac{|w_{i+1} {-} w_i|}{\tau}\right) - C
\ee
with a constant $C>0$ independent of $x$ and $\tau$. 

Coming back to \eqref{dp2}, we note that the boundary source term
$$
\ipo \gamma(x)(u^*_{i+1} - u^*_i)f(w_{i+1} - w_i) \dd s(x)
$$
is bounded by a constant by virtue of Hypothesis~\ref{hy2}, while the boundary term on the left-hand side
$$
\ipo \gamma(x)(u_{i+1} - u_i)f(w_{i+1} - w_i) \dd s(x)
$$
is non-negative by monotonicity of both functions $f$ and $g$. Using \eqref{dp4} we thus get
\begin{align}\nonumber
&\sumim \io |w_{i+1} - w_i|\log\left(1 + \frac{|w_{i+1} {-} w_i|}{\tau}\right)\dd x\\ \label{dp5}
&\qquad + \sumim\io \Big(\kappa(x,s_{i+1})\nabla u_{i+1}-\kappa(x,s_i) \nabla u_i\Big)\cdot \nabla f(w_{i+1} - w_i)\dd x \le C
\end{align}
with a constant $C>0$ independent of $\tau$. We further have
\begin{align*}
&\Big(\kappa(x,s_{i+1})\nabla u_{i+1}-\kappa(x,s_i) \nabla u_i\Big)\cdot \nabla f(w_{i+1} {-} w_i)\\
&\ \, = f'(w_{i+1} {-} w_i)\Bigg(\Big(\big(\kappa(x,s_{i+1}) {-} \kappa(x,s_{i})\big)\nabla u_i + \kappa(x,s_{i+1})\nabla (u_{i+1} {-} u_i)\Big)\times\\
&\qquad \times \Big(\big(g'(u_{i+1}) {-} g'(u_{i})\big)\nabla u_i + g'(u_{i+1})\nabla (u_{i+1} {-} u_i)\Big)\Bigg)\\
&\ \, = f'(w_{i+1} {-} w_i)\Bigg(\!g'(u_{i+1})\kappa(x,s_{i+1})|\nabla (u_{i+1} {-} u_i)|^2 + \big(g'(u_{i+1}) {-} g'(u_{i})\big)\!\big(\kappa(x,s_{i+1}) {-} \kappa(x,s_{i})\big)\!|\nabla u_i|^2\\
&\qquad + \Big(\kappa(x,s_{i+1})\big(g'(u_{i+1}) {-} g'(u_{i})\big) + g'(u_{i+1})\big(\kappa(x,s_{i+1}) {-} \kappa(x,s_{i})\big)\Big)\nabla u_i\cdot \nabla (u_{i+1} {-} u_i)\Bigg).
\end{align*}
The functions $\kappa$ and $g'$ are bounded and Lipschitz continuous, and
$f'(w_{i+1} {-} w_i) = \tau/(\tau + |w_{i+1} {-} w_i|)^2$. Moreover, since $s_i = P[w]_i$ admits a representation similar to \eqref{de3}, by Hypothesis~\ref{hy2}, estimate \eqref{GiL}, and the Lipschitz continuity of the time-discrete play implied by \eqref{de4a} we obtain
$$
|\kappa(x,s_{i+1}) - \kappa(x,s_{i})| \le \bar{\kappa}|s_{i+1} - s_i| \le C|w_{i+1} - w_i| \ \ale,
$$
whereas, from assumption \eqref{hg},
$$
|g'(u_{i+1}) {-} g'(u_{i})| \le \bar{g}(U)|u_{i+1} {-} u_i| \le \bar{g}(U)(g_*(U))^{-1}|w_{i+1} {-} w_i|.
$$
Thus, by employing the Cauchy-Schwarz and Young's inequalities we get the estimate
$$
\Big(\kappa(x,s_{i+1})\nabla u_{i+1}-\kappa(x,s_i) \nabla u_i\Big)\cdot \nabla f(w_{i+1} {-} w_i) \ge \frac{-\tau C |w_{i+1} {-} w_i|^2|\nabla u_i|^2}{(\tau + |w_{i+1} {-} w_i|)^2} \ge -\tau C |\nabla u_i|^2
$$ 
with a constant $C>0$ independent of $\tau$. As a consequence of \eqref{dp5}, \eqref{energy}, and \eqref{hg} we thus have the crucial estimate
\be{dp6}
\sumim \io |u_{i+1} - u_i|\log\left(1 + \frac{|u_{i+1} {-} u_i|}{\tau}\right)\dd x \le C\left(1+ \tau \sumiz \io|\nabla u_i|^2 \dd x\right) \le C
\ee
with a constant $C>0$ independent of $\tau$.

Estimates \eqref{Lam}, \eqref{energy}, and \eqref{dp6} are enough to pass to the limit as $\tau\to 0$, provided they ensure suitable compact embeddings to deal with the nonlinear terms. These do not follow in a straightforward way, due to the low (a bit more than $L^1$) regularity of the time derivative provided by estimate \eqref{dp6}, so we cannot apply the theory of anisotropic Sobolev spaces which was used in \cite{colli}. We therefore devote the next two sections to some auxiliary material on compact anisotropic embeddings involving also Orlicz spaces, with the aim of proving a compactness result (namely, Proposition~\ref{pp1}) suitable for our purposes.


\section{Anisotropic Hilbert space embeddings}\label{hilb}

Consider the spaces $H = L^2(0,T)$, $V = W^{1,2}_0(0,T)$, $V^* = W^{-1,2}(0,T)$ endowed with norms
$$
|v|_H = \left(\int_0^T|v(t)|^2\dd t\right)^{1/2},\ |v|_V = \left(\int_0^T|\dot v(t)|^2\dd t\right)^{1/2},\ |v|_{V^*} = \sup_{|w|_V \le 1}\scal{v,w},
$$
where $\left\langle v,w\right\rangle$ for $v \in V^*$ and $w \in V$ denotes the duality between $V$ and $V^*$. The eigenfunctions $\{\ell_j : j\in \nat\}$ of the operator $-(\dd/\dd t)^2$ with zero Dirichlet boundary conditions are given by the formula
\be{ei1}
\ell_j(t) = \sqrt{\frac2T} \sin\left(\frac{j\pi t}{T}\right)
\ee
and form an orthonormal basis in $H$. More specifically, we have
\be{ei2}
-\ddot \ell_j(t) = \mu_j \ell_j(t), \quad \mu_j = \frac{\pi^2 j^2}{T^2}, \quad j\in\nat.
\ee
A function $v \in H$ admits a representation in terms of the Fourier series
\be{fou1}
v(t) = \sumj v_j \ell_j(t), \quad v_j = \int_0^T v(t)\ell_j(t)\dd t.
\ee
For $v \in V$ we have
\be{ei3}
|v|_H^2 = \sumj |v_j|^2, \quad |v|_V^2 = \sumj \mu_j |v_j|^2, \quad |v|_{V^*} = \sup_{\sumj \mu_j|w_j|^2\le 1} \sumj v_j w_j = \left(\sumj \frac1{\mu_j} |v_j|^2\right)^{1/2}.
\ee
The last identity follows from the fact that the maximum of $\sumj s_j r_j$ over all $r_j$ such that $\sumj r_j^2 \le 1$ is reached for $r_j = s_j/\sqrt{\sumi s_i^2}$.

Similarly, in $L^2(\Omega)$ we define an orthonormal basis $\{e_k : k\in \nato\}$ of eigenvectors of the Laplace operator with homogeneous Neumann boundary conditions, that is,
\be{ee1}
\io \nabla e_k\cdot\nabla \vrt \dd x = \omega_k \io e_k\vrt \dd x \quad \forall \vrt \in W^{1,2}(\Omega),
\ee
and such that $\io e_j e_k\dd x = \delta_{jk}$ for all $j,k \in \nato$. The sequence $\{\omega_k : k\in\nato\}$ is nondecreasing, $0 = \omega_0 <\omega_1< \omega_2 \le ...$, $\lim_{k\to \infty} \omega_k = \infty$. 

Let $\HH = L^2(\Omega\times (0,T))$, $\HN = L^2(\Omega\times (0,T);\real^N)$. Then the system $\{\hejk(x,t) = e_k(x)\ell_j(t): (j,k) \in \nat\times(\nato)\}$ forms an orthonormal basis in $\HH$. Let us denote
\begin{equation}\label{XY}
	\XX \coloneqq \{u \in \HH: \nabla u \in \HN\},\ \quad \YY \coloneqq L^2(\Omega; V^*).
\end{equation}
Putting $u_{jk} = \int_0^T\io u(x,t)\hejk(x,t)\dd x\dd t$ for a generic element $u\in \HH$, we obtain the representation
$$
u(x,t) = \sumk\sumj u_{jk}\hejk(x,t), \quad \|u\|_\HH^2:= \int_0^T\io |u(x,t)|^2\dd x\dd t = \sumk\sumj |u_{jk}|^2 < \infty.
$$
According to the above notation, $u \in \XX$ if and only if
\be{ee2}
\|u\|_{\XX}^2 \coloneqq \int_0^T\io \big(|u|^2 + |\nabla u|^2\big)\dd x\dd t = \sumk\sumj (1+\omega_k) |u_{jk}|^2 < \infty,
\ee
and for each $u \in \HH$ we have
\be{ee3}
\|u\|_\YY^2 := \io |u(x,\cdot)|_{V^*}^2 \dd x = \sumk\sumj \frac1{\mu_j} |u_{jk}|^2.
\ee

\begin{lemma}\label{lem1}
The space $\XX$ is compactly embedded in $\YY$.
\end{lemma}

Lemma~\ref{lem1} is compatible with the general theory of anisotropic Sobolev spaces as in \cite{bin}. The Hilbert case, however, can be treated in an elementary way and we show the details of the proof for the reader's convenience.

\bpf{Proof}
Let $A\subset \XX$ be a bounded set. In agreement with \eqref{ee2} this means that there exists $C_A>0$ such that
\be{ei4}
\sumk\sumj (1+\omega_k) |u_{jk}|^2 \le C_A \quad \forall u \in A.
\ee
Let $\ve > 0$ be arbitrarily given. Compactness will be proved if we check that there exists a finite set $\{w^{(1)}, \dots w^{(n)}\} \subset A$ such that 
\be{ei5}
\forall u \in A \quad \exists i \in \{1, \dots, n\}\quad \|u - w\oi\|_\YY^2 < \ve.
\ee
Notice that, without loss of generality, we can assume that $A$ is closed. Indeed, a finite $\ve$-covering of $A$ generates a $2\ve$-covering of $\bar A$. We find $m \in \nat$ such that
\be{er}
\omega_m > 4C_A/(\mu_1\ve),\quad \mu_m \ge 4C_A/\ve,
\ee
and split $\HH$ into the direct sum $\HH = \HH_m \oplus \HH_m^\perp$, where
$$
\HH_m = \sspan\{\hejk: j,k \le m\},\quad \HH_m^\perp = \sspan\{\hejk: j\ge m+1 \ \mbox{ or }\ k\ge m+1 \}.
$$
Since $A$ is closed, the set $A_m:= A\cap \HH_m$ is a bounded closed subset of the $m\times(m+1)$-dimensional space $\HH_m$. Hence, it is compact in $\HH_m$, and therefore admits a finite subset $\{w^{(1)}, \dots w^{(n)}\} \subset A_m$ such that for every $u \in A_m$ there exists $i \in \{1, \dots, n\}$ such that
\be{ehm}
\|u - w\oi\|_\YY^2 = \sum_{k=0}^m\sum_{j=1}^m \frac1{\mu_j}|u_{jk} - w\oi_{jk}|^2 < \frac\ve{2},
\ee
where $w\oi_{jk} = \int_0^T\io w\oi(x,t)\hejk(x,t)\dd x\dd t$.

Let now $u \in A$ be arbitrary. We have
$$
\|u - w\oi\|_\YY^2 = \sum_{k=0}^m\sum_{j=1}^m \frac1{\mu_j}|u_{jk} - w\oi_{jk}|^2 + \sum_{k=m+1}^\infty\sumj\frac1{\mu_j}|u_{jk}|^2 + \sumk\sum_{j=m+1}^\infty\frac1{\mu_j}|u_{jk}|^2,
$$
where we have by virtue of \eqref{ei4}, \eqref{er}, and \eqref{ehm} that
\begin{align*}
\sum_{k=0}^m\sum_{j=1}^m\frac1{\mu_j}|u_{jk}-w\oi_{jk}|^2 &< \frac\ve{2},\\
\sum_{k=m+1}^\infty\sumj\frac1{\mu_j}|u_{jk}|^2 &\le \frac1{\mu_1(1+\omega_m)}\sumk\sumj (1+\omega_k) |u_{jk}|^2 \le \frac{C_A}{\mu_1\omega_m} < \frac\ve{4},\\
\sumk\sum_{j=m+1}^\infty\frac1{\mu_j}|u_{jk}|^2 &\le \frac1{\mu_m}\sumk\sumj|u_{jk}|^2 \le \frac{C_A}{\mu_m} < \frac\ve{4},
\end{align*}
which we wanted to prove.
\epf


\section{Orlicz spaces}\label{orli}

In this section, we mainly refer to \cite{leo}, see also \cite{ada,rare}. A function $\Phi:[0,\infty) \to [0,\infty)$ is called a strict Young function if it is convex, $\lim_{u\to 0}\Phi(u)/u = 0$, and $\lim_{u\to \infty}\Phi(u)/u = \infty$. Let us first recall a classical and well-known result.

\begin{lemma}\label{ol2}
Let $\vp:[0,\infty) \to [0,\infty)$ be an increasing function such that $\vp(0) = 0$ and $\lim_{r \to \infty} \vp(r) = \infty$, and let $\vp^{-1}$ be its inverse. Then for all $u\ge 0$, $v\ge 0$ we have
\be{o1}
\int_0^u \vp(r)\dd r + \int_0^v \vp^{-1}(s)\dd s = uv + \int_{\vp^{-1}(v)}^u (\vp(r) - v)\dd r \ge uv.
\ee
\end{lemma}

\bpf{Proof}
Assume first that $\vp$ is continuously differentiable. The general case is then obtained by approximation. By substituting $s = \vp(r)$ and integrating by parts we have
$$
\int_0^v \vp^{-1}(s)\dd s = \int_0^{\vp^{-1}(v)} r \vp'(r)\dd r = \vp^{-1}(v)v - \int_0^{\vp^{-1}(v)}\vp(r)\dd r,
$$
and \eqref{o1} follows easily.
\epf

For $\vp$ as in Lemma~\ref{ol2} put $\Phi(u) = \int_0^u \vp(r)\dd r$, $\Phi^*(v) = \int_0^v \vp^{-1}(s)\dd s$. From Lemma~\ref{ol2} it follows that
\be{o2}
\Phi^*(v) = \sup_{u\ge 0}\{uv - \Phi(u)\} = v\vp^{-1}(v) - \Phi(\vp^{-1}(v)).
\ee
Under the hypotheses of Lemma~\ref{ol2}, both $\Phi$ and $\Phi^*$ are strict Young functions, and $\Phi^*$ is called the conjugate of $\Phi$. Clearly, if $\Phi^*$ is the conjugate of $\Phi$, then $\Phi$ is the conjugate of $\Phi^*$. Classical examples of conjugate Young functions are
\begin{align}\label{ox1}
&\Phi(u) = \frac1p u^p,\ \Phi^*(v) = \frac1{q}v^{q}, \quad \frac1p + \frac1{q} = 1,\\ \label{ox2}
&\Phi(u) = (1+u)\log(1+u)-u,\ \Phi^*(v) = \expe^v - v - 1.
\end{align}
In terms of the function $\vp$ as in Lemma~\ref{ol2}, the case \eqref{ox1} corresponds to the choice $\vp(r) = r^{p-1}$, while in \eqref{ox2} we have $\vp(r) = \log(1+r)$. Given two Young functions $\Phi$ and $\hat\Phi$, we say that $\Phi$ is dominant over $\hat \Phi$ if there exist $\alpha,\beta >0$ such that $\Phi(u) \ge \frac{\alpha}{\beta} \hat \Phi(\beta u)$ for all $u>0$. Indeed, in this case the conjugate functions satisfy the inequality $\hat \Phi^*(v) \ge \frac{\beta}{\alpha} \Phi^*(\alpha v)$ for all $v>0$, hence $\hat \Phi^*$ is dominant over $\Phi^*$. We say that $\Phi$ and $\hat \Phi$ are equivalent if $\Phi$ is dominant over $\hat \Phi$ and $\hat \Phi$ is dominant over $\Phi$.

For a strict Young function $\Phi$ and a measurable set $Z \subset \real^d$ for $d \in\nat$ we define the sets
\begin{align*}
L^\Phi(Z) &= \left\{u \in L^1(Z):\ \exists a>0 \ \int_Z \Phi(a |u(z)|)\dd z < \infty\right\},\\
L^\Phi_{\sharp}(Z) &= \left\{u \in L^1(Z):\ \forall a>0 \ \int_Z \Phi(a |u(z)|)\dd z < \infty\right\}.
\end{align*}
The sets $L^\Phi(Z)$ and $L^\Phi_{\sharp}(Z)$ are called the large (small, respectively) Orlicz classes on $Z$ generated by $\Phi$. Clearly, Orlicz classes generated by equivalent Young functions coincide.

In the situation of \eqref{ox1} and $Z = (0,T)$ we have $L^\Phi(0,T)=L^\Phi_{\sharp}(0,T) = L^p(0,T)$ and $L^{\Phi^*}(0,T)=L^{\Phi^*}_{\sharp}(0,T) = L^{q}(0,T)$, while in \eqref{ox2} only $L^\Phi(0,T)=L^\Phi_{\sharp}(0,T)$ holds, while $L^{\Phi^*}(0,T)\ne L^{\Phi^*}_{\sharp}(0,T)$. Indeed, for the function $u(t) = -\log t$ and $T=1$ we have
$$
\int_0^1 \Phi^*(a u(t))\dd t < \infty\ \for\ a<1,\quad \int_0^1 \Phi^*(a u(t))\dd t = \infty\ \for\ a\ge 1.
$$ 

For a general strict Young function $\Phi$ we define on $L^\Phi(Z)$ and $L^\Phi_{\sharp}(Z)$ the Luxemburg norm
\be{o3}
|u|_{\Phi,Z} = \inf\left\{b>0: \int_Z \Phi\left(\frac{|u(z)|}{b}\right)\dd z \le 1\right\}.
\ee
Both $L^\Phi(Z)$ and $L^\Phi_{\sharp}(Z)$ endowed with the Luxemburg norm are Banach spaces. Here again, the Young function $\Phi$ in \eqref{ox1} leads to the usual $L^p$-norm $|u|_{\Phi,Z} = \|u\|_p = \left(\io|u(z)|^p \dd z\right)^{1/p}$ on $L^p(Z)$ (analogously for $\Phi^*$ and $L^{q}$). Moreover, if $u \in L^\Phi(Z)$ and $v \in L^{\Phi^*}(Z)$, then $uv \in L^1(Z)$ and the following counterpart of H\"older's inequality for $L^p$ spaces holds:
\be{hol}
\int_Z u(z)v(z)\dd z \le 2|u|_{\Phi,Z}|v|_{\Phi^*,Z}.
\ee 

\begin{lemma}\label{ol4}
Let $a>0$ be given, let $\Phi$ be a strict Young function, and let $|u|_{\Phi,Z} \ge a$ for some $u \in L^\Phi(Z)$. Then 
\be{o6}
|u|_{\Phi,Z} \le a \int_Z \Phi\left(\frac{|u(z)|}{a}\right)\dd z.
\ee
\end{lemma}

\bpf{Proof}
Note first that the function $v\mapsto \Phi(v)/v$ is nondecreasing in $(0,\infty)$. Indeed, the convexity of $\Phi$ implies that $v\Phi''(v) \ge 0$ for all $v>0$, hence, $v\Phi''(v) + \Phi'(v) \ge \Phi'(v)$, and integrating this inequality we obtain the assertion.

Let now $u \in L^\Phi(Z)$ be such that $b = |u|_{\Phi,Z} \ge a$, and let $M = \{z\in Z: u(z)\ne 0\}$. By definition we have
$$
1 = \int_Z \Phi\left(\frac{|u(z)|}{b}\right)\dd z \le \int_M \frac{|u(z)|}{b}\Phi\left(\frac{|u(z)|}{a}\right)\frac{a}{|u(z)|}\dd z = \frac{a}{b} \int_Z \Phi\left(\frac{|u(z)|}{a}\right)\dd z,
$$
which is precisely \eqref{o6}.
\epf

Note that for every strict Young function $\Phi$ we have (see \cite{ada,leo}) the duality
\be{dua}
L^\Phi(Z) = (L^{\Phi^*}_\sharp(Z))^*.
\ee

We conclude this section by proving a compactness result for the case $Z = (0,T)$.

\begin{lemma}\label{ol3}
Let $\Phi$ be an arbitrary strict Young function. Then the space
$$
W^{1,\Phi}(0,T) = \{u \in L^1(0,T): \dot u \in L^\Phi(0,T)\}
$$
endowed with the norm
$$
|u|_{1,\Phi} = \int_0^T |u(t)|\dd t + |\dot u|_{\Phi,(0,T)}
$$
is compactly embedded in the space $C[0,T]$ of continuous functions on $[0,T]$.
\end{lemma}

\bpf{Proof}
As is customary, we employ the Arzel\`{a}-Ascoli Theorem to prove the compact embedding in the space of continuous functions. Consider a bounded subset $K\subset W^{1,\Phi}(0,T)$. We first prove that $\sup\{|u(t)|: u \in K, \ t\in [0,T]\}< \infty$. Indeed, for every $s,t \in [0,T]$ we have
\be{o4}
|u(t)| \le |u(s)| + \int_0^T |\dot u(r)|\dd r.
\ee
By integrating \eqref{o4} with respect to $s$ from $0$ to $T$ and using \eqref{hol} it follows that
$$
\sup_{t \in [0,T]}|u(t)| \le \frac1T\int_0^T |u(s)|\dd s + 2|\dot u|_{\Phi,(0,T)}|1|_{\Phi^*,(0,T)}
$$
which yields the desired upper bound. It remains to prove the equicontinuity
\be{o5}
\forall \ve > 0\ \exists \delta>0\ \forall u \in K: |t-s|\le \delta \ \Longrightarrow\ |u(t) - u(s)| \le \ve.
\ee
We use again \eqref{o4} and \eqref{hol} to obtain the inequality
$$
|u(t) - u(s)| \le 2|\dot u|_{\Phi, (s,t)}\,\inf\left\{b>0: |t-s|\Phi^*\left(\frac1b\right)\le 1\right\} = \frac{2|\dot u|_{\Phi, (s,t)}}{(\Phi^*)^{-1}(1/|t-s|)},
$$
and \eqref{o5} follows.
\epf


\section{Proof of Theorem~\ref{t1}}\label{proo}

We construct the sequences $\{\hat u\on\}$ and $\{{\bar u}^{(n)}\}$ of approximations by piecewise linear and piecewise constant interpolations of the solutions to the time-discrete problem associated with the division $t_i=i\tau=iT/n$ for $i=0,1,\dots,n$ of the interval $[0,T]$. Namely, for $x \in \Omega$ and $t\in [t_{i-1},t_i)$, $i=1, \dots, n$, we define 
\begin{align}\label{pl}
\hat u\on(x,t) &= u_{i-1}(x) + \frac{t-t_{i-1}}{\tau} (u_i(x) - u_{i-1}(x)),\\
\label{pco}
\bar u\on(x,t) &= u_{i}(x),
\end{align}
continuously extended to $t=T$. Similarly, we consider
\begin{align*}
\hat G\on(x,t) &= G[u]_{i-1}(x) + \frac{t-t_{i-1}}{\tau} (G[u]_i(x) - G[u]_{i-1}(x)),\\
\bar G\on(x,t) &= G[u]_i(x),
\end{align*}
and the boundary data ${\bar u}^{*(n)}(x,t) = u^*(x,t_{i}) = u^*_{i}(x)$ for $x \in \partial\Omega$, again continuously extended to $t=T$. Then \eqref{dis1} is of the form
\be{de0w}
\io \left(\hat G\on_t\vrt + \kappa(x,\bar G\on)\nabla \bar u\on \cdot\nabla\vrt\right)\dd x + \ipo \gamma(x)(\bar u\on - {\bar u}^{*(n)})\vrt \dd s(x) = 0
\ee
for every test function $\vrt \in W^{1,2}(\Omega)$. Our goal is to let $n \to \infty$ (or, equivalently, $\tau \to 0$) in \eqref{de0w}, and obtain in the limit a solution to \eqref{e4}--\eqref{e5}. This will be achieved by exploiting the a~priori estimates derived in the previous sections.

The function $\Phi_{log}(v) = v\log(1+v)$ which appears in Theorem~\ref{t1} is a strict Young function. An easy computation shows that it is equivalent to the function $\Phi$ in \eqref{ox2}. More precisely, comparing the second derivatives of $\Phi$ and $\Phi_{log}$ we get $\Phi(u) \le \Phi_{log}(u) \le 2 \Phi(u)$ for all $u>0$. The space $W^{1,\Phi_{log}}(0,T)$ is compactly embedded in $C[0,T]$ by virtue of Lemma~\ref{ol3}. Furthermore, we have the implication
\be{ep1}
\int_0^T\io |u_t|\log(1+ |u_t|) \dd x\dd t \le C \ \Longrightarrow \ \io |u_t(x,\cdot)|_{\Phi_{log},(0,T)}\dd x \le C+|\Omega|.
\ee
Indeed, we split $\Omega$ into $\Omega_+ \coloneqq \{x\in \Omega: |u_t(x,\cdot)|_{\Phi_{log},(0,T)} \ge 1\}$ and $\Omega_- \coloneqq \Omega\setminus \Omega_+$. We have
$$
\io |u_t(x,\cdot)|_{\Phi_{log},(0,T)}\dd x = \int_{\Omega_+} |u_t(x,\cdot)|_{\Phi_{log},(0,T)}\dd x + \int_{\Omega_-} |u_t(x,\cdot)|_{\Phi_{log,(0,T)}}\dd x
$$ 
and \eqref{ep1} follows from Lemma~\ref{ol4} for $a=1$.

From the estimates \eqref{Lam}, \eqref{energy}, and \eqref{dp6}, and thanks to the implication \eqref{ep1} which holds true for the functions $\hat u\on_t$ defined in \eqref{pl} as indeed
\begin{equation}\label{ep1a}
	\int_0^T\io |\hat u\on_t|\log(1+|\hat u\on_t|) \dd x\dd t = \sum_{i=1}^n\int_{(i-1)\tau}^{i\tau}\io \frac{|u_i - u_{i-1}|}{\tau}\log\left(1+\frac{|u_i - u_{i-1}|}{\tau}\right) \dd x\dd t \le C,
\end{equation}
we conclude that the sequence $\hat u\on$ is bounded in the space
\be{ep2}
\BB:= L^\infty(\Omega\times (0,T))\cap \XX \cap L^1(\Omega; W^{1,\Phi_{log}}(0,T))
\ee
with $\XX$ as in \eqref{XY}. We have in particular
\be{ep2a}
\supess_{(x,t) \in \Omega\times (0,T)} |\hat u\on(x,t)| + \|\hat u\on\|_\XX + \io |\hat u_t\on(x,\cdot)|_{\Phi_{log},(0,T)}\dd x \le R
\ee
with a constant $R>0$ independent of $n$, where the norm $\|\,\cdot\,\|_\XX$ is defined in \eqref{ee2}.

By Lemma~\ref{ol3} we have the compact embedding
\be{comp}
W^{1,\Phi_{log}}(0,T) \hookrightarrow\hookrightarrow C[0,T] \hookrightarrow H \hookrightarrow\hookrightarrow V^*
\ee
which enables us to prove the compactness of the sequence of approximate solutions following the idea of \cite{lions}, where we interchange the role of time and space variables. The result reads as follows.

\begin{proposition}\label{pp1}
Let $K_R \subset \BB$ be the set of all $u \in \BB$ satisfying the condition \eqref{ep2a}. Then $K_R$ is compact in the space $L^1(\Omega; C[0,T])$.
\end{proposition}

\bpf{Proof}
The statement will be proved if for every $\ve>0$ we find a finite set $\{w^{(1)}, \dots, w^{(k)}\} \subset K_R$ with the property
\be{dens}
\forall u\in K_R\ \exists j \in \{1, \dots, k\}: \io \max_{t \in [0,T]}|u(x,t) - w\oj(x,t)|\dd x < \ve. 
\ee
As a first step towards the proof of \eqref{dens}, we prove the following statement.
\be{ep3}
\forall \eta > 0 \ \exists L_\eta > 0 \ \forall u \in W^{1,\Phi_{log}}(0,T): \max_{t \in [0,T]}|u(t)| \le \eta |u|_{1,\Phi_{log}} + L_\eta |u|_{V^*}.
\ee
We proceed by contradiction. Assume that
\be{ep4}
\exists \eta>0 \ \forall m\in \nat \ \exists u_m \in W^{1,\Phi_{log}}(0,T): \max_{t \in [0,T]}|u_m(t)| > \eta |u_m|_{1,\Phi_{log}} + m |u_m|_{V^*}.
\ee
We can assume that $|u_m|_{1,\Phi_{log}} = 1$ for all $m \in \nat$. By Lemma~\ref{ol3}, $W^{1,\Phi_{log}}(0,T)$ is compactly embedded in $C[0,T]$. There exists therefore $u \in C[0,T]$ such that a subsequence $\{u_{m_i}\}$ of $\{u_m\}$ converges uniformly to $u$. Since $C[0,T]$ is continuously embedded in $V^*$, $u_{m_{i}}$ converge to $u$ in $V^*$. Hence, by \eqref{ep4}, $u = 0$, which contradicts the assumption that $|u_m|_{1,\Phi_{log}} = 1$ for all $m \in \nat$. Hence, \eqref{ep3} is proved.

Let us pass to the proof of \eqref{dens}. Let $\ve > 0$ be given.
In \eqref{ep3} we choose $\eta = h\ve$ for $h>0$ which will be specified later. By Lemma~\ref{lem1}, the space $\XX$ is compactly embedded in $\YY$. Hence, the closure $\bar K_R^{\XX}$ of $K_R$ in $\XX$ is compact in $\YY$, and we can find a finite set $\{w^{(1)}, \dots, w^{(k)}\} \subset \bar K_R^{\XX}$ such that for every $u \in \bar K_R^{\XX}$ there exists $j \in \{1,\dots, k\}$ such that $\|u - w\oj\|_\YY < \eta/L_\eta$. Since $K_R$ is dense in $\bar K_R^{\XX}$, we can assume that $w\oj \in K_R$ for all $j \in \{1,\dots, k\}$. For each $u \in K_R$ we integrate \eqref{ep3} for $u - w\oj$ instead of $u$ over $x \in \Omega$, and using \eqref{ep2a} we obtain the inequality
\begin{align*}
&\io\max_{t \in [0,T]} |u(x,t) - w\oj(x,t)|\dd x\\
 &\qquad \le \eta\left(\int_0^T\io(|u|+|w\oj|)(x,t)\dd x\dd t+ \io (|u_t(x,\cdot)|_{\Phi_{log}} + |w\oj_t(x,\cdot)|_{\Phi_{log}}\big)\dd x\right)\\
&\hspace{12mm} + L_\eta \io |u(x,\cdot) - w\oj(x,\cdot)|_{V^*}\dd x\\
&\qquad \le 2\big(\sqrt{|\Omega|T}+1\big)R\eta + \sqrt{|\Omega|}L_\eta \|u-w\oj\|_\YY \\
&\qquad \le \left(2\big(\sqrt{|\Omega|T}+1\big)R + \sqrt{|\Omega|}\right)\eta.
\end{align*}
Choosing $h$ such that $h\Big(2\big(\sqrt{|\Omega|T}+1\big)R + \sqrt{|\Omega|}\Big)= 1$, we conclude the proof of Proposition~\ref{pp1}.
\epf

From Proposition~\ref{pp1} and \eqref{ep2a} it follows that the sequence $\hat u\on$ of approximate solutions defined in \eqref{pl} is compact in the space $L^1(\Omega; C[0,T])$. Furthermore, the time derivatives $\hat u\on_t$ are bounded in $L^{\Phi_{log}}(\Omega\times (0,T))$, see \eqref{ep1a}. Hence, a subsequence of $\hat u\on$ converges strongly in $L^1(\Omega; C[0,T])$ to $u \in \BB$ and, by \eqref{dua}, $\hat u\on_t \to u_t$ weakly* in $L^{\Phi_{log}}(\Omega\times (0,T))$. By \cite[formula~(35)]{colli}, we deduce that an estimate analogous to \eqref{ep1a} holds for $\hat G\on_t$ as well, so there exists $G^* \in L^{\Phi_{log}}(\Omega\times (0,T))$ such that $\hat G\on_t \to G^*$ weakly* in $L^{\Phi_{log}}(\Omega\times (0,T))$. The Lipschitz continuity of $G$ in $L^1(\Omega;C[0,T])$ stated in Proposition~\ref{pc1} implies that $G[\hat u\on] \to G[u]$ in $L^1(\Omega;C[0,T])$.

The main issue is the convergence proof of the approximate hysteresis terms in Eq.~\eqref{de0w}. This is obtained by exploiting a result about Lipschitz continuity of $G$ with respect to the sup-norm in the space of regulated functions of time (that is, functions which admit finite left and right limits at every point $t$) in \cite[Proposition~2.3]{jana}. It states that if $v,w:[0,T] \to \real$ are regulated functions, the initial memory curves $\lambda$ for $G[v]$ and $G[w]$ coincide, and $[a,b] \subset [0,T]$ is an arbitrary interval, then there exists a constant $C_G>0$ independent of $v,w,a,b$ such that
$$
\sup_{t\in [a,b]}|G[v](t)- G[w](t)| \le C_G \sup_{t\in [a,b]}|v(t)- w(t)|.
$$
Since in our case the initial memory curve $\lambda$ from Definition~\ref{dpr} is the same for $\hat{G}\on$ and for $G[\hat u\on]$, we have for all $t \in [0,T]$ and a.\,e.\ $x \in \Omega$ that
$$
|\hat G\on(x,t) {-} G[\hat u\on](x,t)| \le C_G \max_{i=1, \dots, n} |u_i(x) {-} u_{i-1}(x)|.
$$
Since $\Phi_{log}(Cv) \le C\Phi_{log}(v)$ for $C\ge 1$, we can assume that $C_G \ge 1$ and get that
\begin{equation}\label{difg}
\begin{aligned}
&\io\Phi_{log}\bigg(\sup_{t \in (0,T)}|\hat G\on {-} G[\hat u\on]|(x,t)\bigg)\dd x \le \io\Phi_{log}\bigg(C_G\max_{i=1, \dots, n} |u_i(x) {-} u_{i-1}(x)|\bigg)\dd x\\
&\qquad \le C_G \sumin\io\Phi_{log}(|u_i(x) {-} u_{i-1}(x)|)\dd x.
\end{aligned}
\end{equation}
By \eqref{Lam} we have $0 \le |u_i(x) - u_{i-1}(x)| \le 2U$ for all $i=1, \dots, n$ and a.\,e.\ $x \in \Omega$. To estimate the right-hand side of \eqref{difg}, we use the following easy Lemma.

\begin{lemma}\label{lphi}
	There exists a function $\alpha:(0,1)\to (0,1)$ such that $\lim_{\tau \to 0} \alpha(\tau) = 0$ and for all $v \in (0,2U)$ and $\tau \in (0,1)$ we have
	$$
	\frac{\Phi_{log}(v)}{\tau\Phi_{log}\left(\frac{v}{\tau}\right)} = \frac{\log(1+v)}{\log\left(1+\frac{v}{\tau}\right)} \le \alpha(\tau).
	$$
\end{lemma} 

\bpf{Proof}
For $v \in (\tau^{1/2}, 2U]$ we have $\log\left(1+\frac{v}{\tau}\right) \ge \log(1+\tau^{-1/2})$, $\log (1+v) \le \log(1+2U)$, and
\be{al1}
\frac{\log(1+v)}{\log\left(1+\frac{v}{\tau}\right)} \le \frac{\log(1+2U)}{\log\left(1+\tau^{-1/2})\right)} =: \alpha_1(\tau).
\ee
For $v \in (\tau, \tau^{1/2}]$ we have $\log\left(1+\frac{v}{\tau}\right) \ge \log 2$, $\log (1+v) \le \tau^{1/2}$, and
\be{al2}
\frac{\log(1+v)}{\log\left(1+\frac{v}{\tau}\right)} \le \frac{\tau^{1/2}}{\log 2} =: \alpha_2(\tau).
\ee
We now continue for $v \in (\tau^{(k+1)/2}, \tau^{k/2}]$ for $k=2,3, \dots$. Then
$\log\left(1+\frac{v}{\tau}\right) \ge \frac{v}{2\tau} \ge \frac12\tau^{(k-1)/2}$, $\log (1+v) \le \tau^{k/2}$, and
\be{alk}
\frac{\log(1+v)}{\log\left(1+\frac{v}{\tau}\right)} \le 2\tau^{1/2}=: \alpha_3(\tau).
\ee
It suffices to put $\alpha(\tau) = \max\{\alpha_1(\tau), \alpha_2(\tau), \alpha_3(\tau)\}$ to complete the proof.
\epf

Using Lemma~\ref{lphi} and inequality \eqref{dp6}, we estimate the right-hand side of \eqref{difg} as
\begin{equation}\label{difg1}
	\sumin\io\Phi_{log}\big(|u_i(x) {-} u_{i-1}(x)|\big)\dd x \le \alpha(\tau) \sumin \io \tau \Phi_{log}\bigg(\frac{|u_i(x) {-} u_{i-1}(x)|}{\tau}\bigg)\dd x \le C\alpha(\tau),
\end{equation}
which implies
\be{difga}
\io\Phi_{log}\bigg(\sup_{t \in (0,T)}|\hat G\on {-} G[\hat u\on]|(x,t)\bigg)\dd x \le C\alpha(\tau)
\ee
and similarly
\be{difgb}
\io\Phi_{log}\bigg(\sup_{t \in (0,T)}|\bar G\on {-} G[\hat u\on]|(x,t)\bigg)\dd x \le C\alpha(\tau),
\ee
with a constant $C>0$ independent of $\tau$.

Recall that $G[\hat u\on]$ converge strongly to $G[u]$ in $L^1(\Omega; C[0,T])$. Both $\hat G\on$ and $\bar G\on$ are bounded by virtue of \eqref{GiL}. Hence, by \eqref{difga}--\eqref{difgb} and the Lebesgue dominated convergence theorem, they converge to $G[u]$ strongly in $L^p(\Omega\times (0,T))$ for all $1\le p < \infty$.

As a last step, we pass to the limit in the diffusion and in the boundary terms. Note that, to this purpose, we also need to prove the convergence of $\bar{u}\on$ to the same limit $u$. We first notice that
for a.\,e.\ $x\in\Omega$ and all $t\in [0,T]$ we have
$$
|\hat u\on(x,t) - \bar u\on(x,t)| \le \max_{i=1, \dots, n} |u_i(x) - u_{i-1}(x)|,
$$
hence, 
\be{difp}
\sup_{t\in (0,T)}\Phi_{log}\big(|\hat u\on(x,t) - \bar u\on(x,t)|\big) \le \Phi_{log}\bigg(\max_{i=1, \dots, n} |u_i(x) - u_{i-1}(x)|\bigg) \le \sumin \Phi_{log}\big(|u_i(x) - u_{i-1}(x)|\big).
\ee
Using inequality \eqref{difg1} we estimate
\be{difphi}
\io \sup_{t\in (0,T)}\Phi_{log}\big(|\hat u\on(x,t) - \bar u\on(x,t)|\big)\dd x \le C\alpha(\tau),
\ee
with a constant $C>0$ independent of $\tau$. Since $\hat{u}\on$ converge strongly to $u$ in $L^1(\Omega;C[0,T])$ and $\bar{u}\on$ is uniformly bounded by virtue of \eqref{Lam}, by \eqref{difphi} and again the Lebesgue dominated convergence theorem we see that $\bar{u}\on \to u$ strongly in $L^p(\Omega\times (0,T))$ for all $1\le p < \infty$. Hence, by \eqref{energy}, $\nabla\bar{u}\on \to \nabla u$ weakly in $L^2(\Omega\times (0,T);\real^N)$. Moreover, since we have proved that $\bar{G}\on \to G[u]$ strongly in $L^p(\Omega\times (0,T))$ for all $1\le p < \infty$ and $\kappa$ is bounded and (Lipschitz) continuous by Hypothesis~\ref{hy2}, by the Lebesgue dominated convergence theorem we also obtain that $\kappa(\cdot,\hat{G}\on) \to \kappa(\cdot,G[u])$ strongly in $L^p(\Omega\times (0,T))$ for all $1\le p < \infty$. This is enough to pass to the (weak) limit in the diffusion term $\kappa(x,\bar G\on)\nabla \bar u\on$. As for the boundary term, we apply the interpolation inequality
$$
\ipo |v|^2 \dd s(x) \le C\left(\io |v|^2 \dd x + \io |v||\nabla v| \dd x\right) \quad \mbox{for } v \in W^{1,2}(\Omega)
$$
to $v = \bar{u}\on-u$, and conclude that $\bar{u}\on \to u$ strongly in $L^2(\partial\Omega\times (0,T))$. Finally, by Hypothesis~\ref{hy2} the boundary data ${\bar u}^{*(n)}$ are uniformly bounded, hence they converge strongly to $u^*$ (at least) in $L^2(\partial\Omega\times (0,T))$.

To pass to the limit in \eqref{de0w}, we choose any smooth test function $\sigma(t)$ with compact support in $(0,T)$ (note that such functions form a dense	subset of $L^2(0,T)$) and integrate by parts to get
\be{weakn}
\int_0^T\!\!\io \big(-\dot\sigma(t)\hat G\on\vrt + \sigma(t)\kappa(x,\bar G\on)\nabla \bar u\on \cdot\nabla\vrt\big)\dd x\dd t {+} \int_0^T\!\!\ipo \sigma(t)\gamma(x)(\bar u\on {-} {\bar u}^{*(n)})\vrt \dd s(x)\dd t = 0.
\ee
The limit as $n \to \infty$ has the form
\be{weak}
\int_0^T\io \big(-\dot\sigma(t) G[u]\vrt + \sigma(t)\kappa(x,G[u])\nabla u\cdot\nabla\vrt\big)\dd x\dd t + \int_0^T\ipo \sigma(t)\gamma(x)(u - u^*)\vrt \dd s(x)\dd t = 0.
\ee
This also proves that $G^* = G[u]_t \in L^{\Phi_{log}}(\Omega\times (0,T)) \subset L^1(\Omega\times (0,T))$. Thus, by restricting ourselves to test functions $\vrt \in W^{1,2}(\Omega) \cap L^{\infty}(\Omega)$, we have completed the proof of Theorem~\ref{t1}.


\end{document}